\documentclass[11pt]{amsart}
\usepackage{amsmath,amsfonts,euscript,amscd,amsthm,amssymb,upref,graphicx, }

\theoremstyle{plain}
\newtheorem{theorem}{Theorem}

\newtheorem{definition}{Definition}
\newtheorem{proposition}{Proposition}
\newtheorem{lemma}{Lemma}

\swapnumbers

\theoremstyle{definition}
\newtheorem{example}[subsection]{Example}
\newtheorem{remark}[subsection]{Remark}

\newtheorem{nothing*}[subsection]{}

\newcommand{\rien}[1]{}

\newcommand{\dL}{ \operatorname{\mathcal L_0}}
\newcommand{\lL}{ \operatorname{\mathcal L}}
\newcommand{\jL}{ \operatorname{\mathcal L_J}}

\newcommand{\C}{\ensuremath{\mathbb{C}}}

\newcommand{\Ker}{{\rm Ker} \,}

\newcommand{\nn}{\nonumber}

\renewcommand{\epsilon}{\varepsilon}
\renewcommand{\phi}{\varphi}

\begin{document}
\title[Lie algebras of zero divergence vector fields on complex affine algebraic varieties] {Lie algebras of zero divergence vector fields on complex affine algebraic varieties}

\author{Fabrizio Donzelli}
\maketitle \vfuzz=2pt

\begin{abstract}
For a smooth manifold $X$ equipped with a volume form, let
$\dL$ be the Lie algebra of volume preserving smooth vector fields on $X$.
A. Lichnerowicz proved that the abelianization of $\dL$ is a finite-dimensional vector space, and that its dimension
depends only on the topology of $X$.
 In this paper we provide analogous results for some classical examples of non-singular complex affine algebraic varieties that admit a nowhere-zero algebraic form of top degree
(which plays the role of a volume form).
\end{abstract}

\section{Introduction}

We will denote by $X$ (or any other capital letter) a non-singular affine  algebraic variety defined over $\C$, of dimension $\dim X = n$.  The Lie algebra of vector fields on $X$,
which is the same as the Lie algebra of derivations on the ring $\C [X]$ of the regular functions on $X$, will be denoted by  $V_X$.
Let $\Omega^k_X$ be the space of  algebraic $k$-forms, with differential
 $d^k :\Omega^k_X\rightarrow \Omega^{k+1}_X$. The subspaces of closed $k$-forms will be denoted
  by $\Omega^k_{X,cl}$. Since in most of the circumstances there will not be confusion about the variety, 
we will often drop the symbol $X$ from the notation.
For a $k$-form $\eta$ on $X$, the Lie derivative $L_v\eta$ with respect to a vector field $v$ can be  expressed by the Cartan's magic formula

\begin{align}\label{Cartan}
L_v\eta = d i_v \eta+i_v d\eta ,
\end{align}

where $i_v$ is the contraction by $v$. If $X$ has trivial canonical bundle $K_X$, then the $\C [X]$-module
of the global sections of $K_X$ is generated by a nowhere vanishing $n$-form $\omega$, which we will call a volume form on $X$. The reason for the terminology comes from observing that a $n$-dimensional non-singular affine variety admits a natural structure of smooth manifold of real dimension $2n$, and $\omega\wedge\overline{\omega}$ is a smooth volume form on it.
We will assume from now on that $X$ has trivial canonical bundle and that we have fixed a volume form $\omega$.
 Since a $n$-form $\eta$ is uniquely written as $\eta =p\omega$, for some $p\in \C [X]$, it follows that 
the choice of $\omega$ is unique up to multiplication by a unit in $\C [X]$.
 The divergence ${\rm div}_\omega  v$ of a vector field $v$ with respect to $\omega$, which  measures the variation of volume under the flow obtained
 by integrating $v$, is defined by the equation

\begin{align}\label{divergence}
 L_v\omega = ({\rm div}_\omega  v)\omega.
\end{align}

We will denote by $\dL (X)$, or more often by $\dL$, the Lie algebra of zero divergence
 (that is volume preserving) algebraic vector field on $X$.

\begin{lemma}\label{div0}
Let $X$ be a $n$-dimensional non-singular affine algebraic variety equipped with a volume form $\omega$. Then the map
$\mu : V\rightarrow \Omega^{n-1}$ defined by $\mu (v)= i_v\omega$ is an isomorphism of
$\C [X]$-modules. Moreover, $v\in \dL$ if and only if $\mu (v)$ is a closed form.

\begin{proof}
The first statement is a simple fact of linear algebra, the second one follows equations (\ref{Cartan}) and (\ref{divergence})
(more details can be found in 	\cite{KK1}).
\end{proof}

\end{lemma}

\begin{lemma}\label{span}
Let $X$ be a $n$-dimensional non-singular affine algebraic variety  equipped with a volume form $\omega$.
  Let  $\{ v_1,v_2, ... , v_m\}\subset \dL$ be a collection of zero divergence vector fields
such that $\{ v_{1x},v_{2x}, ... , v_{mx}\}$ span the tangent space $T_xX$ for all $x\in X$. Then
$\dL\C[X] =  v_1\C[X]+...+v_m\C [X]$. Moreover, a $n$-form $q\omega $ is exact if and only if
$q\in \dL\C [X]$.

\begin{proof}
Let $\omega_k=i_{v_k}\omega$. Lemma \ref{div0} implies that $\{ \omega_{1x}, ... ,\omega_{mx}\}$ span the fiber $\Omega^{n-1}_x$ for all $x\in X$: if follows from 
\cite{Ha} (Excercise 10, Chapter 2) that 
 the collection $\{ \omega_1, ... ,\omega_m\}$ generates the
$\C [X]$-module $\Omega^{n-1}$.
If $v\in \dL$ and $p\in \C [X]$, a simple calculation shows  that $v(p)\omega = d (pi_v\omega)$.
 Hence we can write
$v(p)\omega = d(p_1\omega_1 + ... p_m\omega_m)= (v_1(p_1)+... + v_m(p_m))\omega$, for some $p_k\in \C [X]$.
It follows that $v(p) = v_1(p_1)+... + v_m(p_m)$, which proves the first assertion. The same argument can be used to prove the second assertion.

\end{proof}

\end{lemma}

The following lemma follow from a straightforward calculation (see \cite{KK1}).
\begin{lemma}\label{comm}
Let $\omega$ be a volume form on the affine algebraic variety $X$. Then for $v,w\in \dL$ the following identity holds:
$$
i_{[v,w]}\omega=di_vi_w\omega
$$

\end{lemma}

We borrow the notion of flexibility, introduced in \cite{AFKKZ}, to highlight Lie algebras of vector fields satisfying the condition of Lemma \ref{span}.

\begin{definition}
A Lie subalgebra $\mathcal G$ of  $V$ is said to be flexible if $\mathcal G_x =T_xX $ for all $x\in X$.
\end{definition}

The following example of a flexible subalgebra will play a crucial role in this paper.

\begin{proposition}\label{jL}
Let $\mu$ be the isomorphism defined in Lemma \ref{div0}. Then  $\jL = \mu^{-1}(d\Omega^{n-2}) $ is a flexible Lie subalgebra of $\dL$. 
In particular, $\dL$ is flexible.

\begin{proof}
Lemma \ref{comm} implies that $\jL$ is a Lie subalgebra of $\dL$. 
In order to show that $\jL$ (and therefore $\dL$) is flexible, embed $X$ as a closed subvariety of $\C^N$ with affine coordinates $(x_1,..., x_N)$.
The collection of exact $(N-1)$ forms $dx_{j_1}\wedge ... \wedge dx_{j_{N-1}}$, for all possible subsets
$\{ j_1,... , j_{N-1} \}$ 
 spans $\Lambda^{N-1}(T_x^*X)$ for all $x\in X$. Hence by Lemma \ref{div0} the collection of vector fields
$v_{j_1,... , j_{N-1}}= \mu^{-1}(dx_{j_1}\wedge ...\wedge dx_{j_{N-1}})$ spans $T_xX$ for all $x\in X$.

\end{proof}

\end{proposition}

\begin{remark}
A vector field  $v\in \jL$ such that $\mu (v) =dp_1\wedge ... \wedge dp_{n-1}$, for $n=\dim X$ is
called a Jacobian vector field (this definition is equivalent the original one given for example in \cite{Fre}).
It is not difficult to show that $\jL$ is generated as a vector space by the set of jacobian vector fields.
Makar-Limanov Theorem estabilishes a connection between jacobian vector fields and the locally nilpotent ones, and perhaps Makar-Limanov's Theorem could explain the phenomenon discussed in Example \ref{matias}
\end{remark}

One of the key ingredients of the paper is a remarkable result of Grothendieck.

\begin{theorem}\label{grot}(\cite{G}, Grothendieck)
Let $X$ be a non-singular affine  algebraic variety defined over $\C$.
Then the de Rham cohomology ring of $X$ is isomorphic to the cohomology ring
of the complex of algebraic differential forms on $X$.
\end{theorem}

Here are some crucial applications of Grothendieck's Theorem.

\begin{proposition}\label{exact}
Let $X$ be a non-singular affine algebraic variety with volume form $\omega$.  If $\mathcal G$ is a flexible
subalgebra of $\dL$, then:
\begin{enumerate}

\item $\dL \C[X]= \mathcal G \C[X]$
\item  $\omega $ is exact if and only if $1\in \mathcal G \C [X]$.
\item $\C[X]/\dL\C [X]= \C [X]/\mathcal G \C[X]\cong {\rm H^{n}}(X,\C)$

\end{enumerate}

\begin{proof}
The first two parts follow from Lemma \ref{span}, while the 
third one is a consequence of Lemma \ref{span} and Theorem
\ref{grot}.

\end{proof}

\end{proposition}

For a Lie algebra $\mathcal G$, we denote by $[\mathcal G , \mathcal G]$ its derived subalgebra,
and the quotient $\mathcal G / [\mathcal G , \mathcal G]$ is called the abelianization of $\mathcal G$.
Lemma \ref{comm} implies that  $[\dL , \dL ] $ is a Lie subalgebra of $\jL$.
Grothendiecks's Theorem and the definition $\jL$ from Proposition \ref{jL} imply the following lemma.

\begin{lemma}\label{hn-1}
Let $X$ be a non-singular affine algebraic variety equipped with a volume form. Then
$\dL /\jL \cong {\rm H}^{n-1}(X, \C)$
\end{lemma}

\subsection{Real differentiable manifolds: Lichnerowicz's work}

Let $M$ be a real orientable smooth manifold (without boundary) of dimension $n$.
It is a standard fact that $M$ admits a nowhere zero smooth $n$-form $\omega$.
The Lie alebras $\dL$ and $\jL$, of respectively zero divergence and jacobian smooth vector fields, 
are defined exactly in the same way as in the previous section.

\begin{theorem}\label{lic}(\cite{L}, Lichnerowicz)
Let $M$ a real smooth orientable manifold. Then $\jL = [\dL , \dL]$.
\end{theorem}

In a survey paper \cite{R} it is claimed that the statement of Theorem  \ref{lic}
is purely algebraic and does not require any smoothness assumption. However, the proof
of Lichnerowicz is based on non-trivial local calculations involving
smooth forms with compact support, and it is not clear how to extend his techniques in the contest
of affine algebraic geometry. One would replace the use of the partition of unity, which is used by Lichnerowicz, by the Hilbert Nullstellensatz, but it is not clear how to do so. 
In this paper we will explore to what extent the theorem of Lichnerowicz
is valid in the contest of affine algebraic geometry. On one hand we will show that Theorem \ref{lic} is true for a large class of affine varieties, on the other hand we will provide a class of affine surfaces for which the Lie algebras $\jL$ and $[\dL , \dL]$ do not coincide.
Unfortunately we can not prove a general statement for affine varieties. However, the examples presented
and the heavy usage of Grothendieck's theorem suggest that $\jL/[\dL , \dL ]$  should be finite-dimensional and that its dimension should depend only on the de Rham cohomology of the affine variety.

\subsection*{Acknowledgments}

The author is very grateful to:
Daniel Daigle, for supporting the postdoctoral fellowship at University of Ottawa and for the interesting conversation that lead to the writing of this paper;
Shulim Kaliman for reading the paper before submission and for important comments;
Pavel Etingof for answering a question on his paper \cite{ES}, mentioned in Remark \ref{Etingof};
Alexander Premet for suggesting the proof of Proposition \ref{semisimple}.

\section{Main results}
In this section we will compute the dimension of $\jL / [ \dL , \dL ]$. Since   $ \dim\dL/\jL =h^{n-1}(X)$ by Lemma \ref{hn-1},
this will give us also the dimension of $\dL / [\dL , \dL ]$.


\begin{proposition}\label{comm}
For the isomorphism $\mu :\dL\rightarrow \Omega^{n-1}_{cl}$ defined in  Lemma \ref{div0},
$\mu ([\dL , \dL ])= \dL\Omega^{n-1}_{cl}$ and
$\jL/[\dL , \dL ]\cong d\Omega^{n-2}/\dL \Omega^{n-1}_{cl}$.

\begin{proof}

Let $w\in\dL $: then $i_w\omega$ is a closed $(n-1)$-form, and the first statement follows from Lemma \ref{Cartan} and Lemma \ref{comm}.
 The second statement follows therefore from the definition of $\jL$.

\end{proof}

\end{proposition}

Cartan's magic formula implies that $L_vd=dL_v$ for any vector field $v$, and therefore 
$d\dL \Omega^{n-2}=\dL d\Omega^{n-2}$ . It follows that the next proposition (which is a consequence of Proposition \ref{comm}) is well stated, since $d\dL \Omega^{n-2}\subset \dL \Omega^{n-1}_{cl}$.

\begin{proposition}\label{dim}
Let $d: \Omega^{n-2}/\dL \Omega^{n-2}\rightarrow d\Omega^{n-2}/d\dL \Omega^{n-2}$ the map
induced by exterior differentiation, 
$\pi: d\Omega^{n-2}/d\dL \Omega^{n-2}\rightarrow d\Omega^{n-2}/\dL \Omega^{n-1}_{cl}$ be the quotient map, and 
$D=\pi\circ d$. If $\Omega^{n-2}/\dL \Omega^{n-2}$ is finite-dimensional then
 $\dim \jL/[\dL , \dL ]=\dim \Omega^{n-2}/\dL \Omega^{n-2}-\dim Ker D $
\end{proposition}

\begin{remark}\label{Etingof}

 Proposition \ref{dim} gives a complete result only in the case of affine algebraic surfaces (Theorem \ref{surface}), 
since for $n=2$  by Proposition \ref{exact} we have that  $\Omega^{n-2}/\dL \Omega^{n-2}= \C [X]/\dL \C [X]\cong {\rm H}^2(X, \C)$
is finite-dimensional. 
It is not clear how to compute the dimension of $\Omega^{n-2}/\dL \Omega^{n-2}$ for $n\geq 3$. A possible method could be found
in \cite{ES}.
\end{remark}

Let $\psi : \Lambda^2 (V)\rightarrow \Omega^{n-2}$ be the map of $\C[X]$-modules
that is defined for indecomposable elements by

\begin{align}
\psi (v\wedge w) =i_vi_w\omega .
\end{align}

A calculation in local coordinates shows that
$\psi$ is an isomorphism. Let $R$ be the image of $\Lambda^2(\dL )$ under $\psi$.
We are going to study the following commutative diagram.

\begin{align}\label{keycd}
\begin{CD}
\Lambda^2(\dL) @>\psi >> R @ > d_{|_R} >> dR =\dL\Omega^{n-1}_{cl}\\
\bigcap @. \bigcap @. \bigcap\\
\Lambda^2(V) @> \psi >> \Omega^{n-2} @> d >> d\Omega^{n-2}
\end{CD}
\end{align}

\vspace{0.5cm}

The commutative diagram allows us to compute the dimension of $\Ker D$ for $n=2$
(where $D$ is defined in Proposition \ref{dim}). It is not clear how to use it in general, nonetheless it will be useful in the case of linear algebraic groups.

\begin{theorem}\label{surface}
Let $X$ be a non-singular affine algebraic surface with a volume form $\omega$. Then

\begin{align}
\nn \dim  \jL/[\dL , \dL ] = h^2(X)
-\dim {\rm span} \{\C [\omega ] , \Lambda^2 ({\rm H}^1 (X, \C))\},
\end{align}

where 
$\C[\omega]\subset {\rm H}^{2}(X,\C) $ is the  space spanned by the cohomology class of $\omega$
and $\Lambda^2 ({\rm H}^1 (X, \C))$ is the subspace  of $H^2(X, \C)$ generated by wedge products of closed 1-forms.

\begin{proof}

Referring to Proposition \ref{dim} we will compute
$\dim\Ker D=\dim \Ker d+\dim \Ker \pi$ for  $n=2$.
For the differential map, we have that  $d: \C[X]/\dL \C [X] \rightarrow d\C [X]/\dL d\C [X]$ and that 
$\dim \Ker d \leq 1$. The second part of Proposition \ref{exact} implies that  $\dim \Ker d = 1$
if and only if the volume form $\omega$ non-exact:  hence $\dim\Ker d=\dim\C [\omega]$.
As for the projection map, in this case we have that $\pi : d\C[X]/\dL d\C [X]\rightarrow  d \C [X]/\dL \Omega^1_{cl}$ and that
$\Ker \pi \cong \dL \Omega^1_{cl}/\dL d\C [X]$. Referring to the first row of 
the commutative diagram \ref{keycd}, $\dL \Omega^1_{cl}/\dL d\C [X]\cong dR/d\dL \C [X]$:
hence $\Ker \pi \oplus (\C\cap R)/ (\C\cap\dL \C [X]) \cong R/\dL \C [X]$.
The proof of the theorem will follow if we show that 
(a) $ R/\dL \C [X] \cong \Lambda^2 (H^1 (X, \C))$ and (b) $(\C\cap R)/ (\C\cap\dL \C [X]) \cong \Lambda^2 (H^1(X, \C))\cap \C [\omega]$.

First we claim that the isomorphism $\psi $ of the commutative diagram \ref{keycd}  restricts to a vector space isomorphism from $\Lambda^2 (\jL )$ to $\dL\C [X]$. In fact,
 note first that if $v,w \in \jL$, then there are functions $f,g\in \C [X]$ such that $i_v\omega =df $ and $i_w\omega =dg$ and therefore  $i_vi_w\omega = v(g)\in \jL \C[X]$.
 Conversely, suppose that $g\in \jL\C[X] $.  Without loss of generality we can assume that $g=v(f)$ and that
$v$ corresponds to an exact form  $i_v\omega $. Moreover there exists a $w\in\jL$ such that $i_w\omega =df$.
Then it is clear that $v\wedge w \in\Lambda^2( \jL)$ and that $\psi (v\wedge w)= v(f)=g$.
Since $\jL \C [X]=\dL \C [X]$ by part (1) of Proposition \ref{exact}, the claim is proved.

This first step allows, in the case $n=2$, to extend the commutative diagram to

\begin{align}\label{keycd2}
\begin{CD}
\Lambda^2 ( \jL )@>\psi >> \dL\C [X] @>d>> d\dL \C [X]\\
\bigcap @. \bigcap @. \bigcap\\
\Lambda^2(\dL) @>\psi >> R @ > d_{|_R} >> dR =\dL\Omega^{1}_{cl}\\
\bigcap @. \bigcap @. \bigcap\\
\Lambda^2(V) @> \psi >> \C [X] @> d >> d\C [X]
\end{CD}
\end{align}

\vspace{0.5cm}

Next we construct a degree-preserving isomorphism between the exterior algebras
$\Omega_X^\bullet$ and $ \Lambda^\bullet(V)$, as follows.
Define $\phi= \phi_0\oplus \phi_1 \oplus \phi_2$, where, on indecomposable elements and then extended by linearity, the three maps $\phi_k$ are defined in the following way:
for $\phi_0:\C[X]\rightarrow \C[X]$, let $\phi_0(g)=g$; for $\phi_1:V\rightarrow \Omega^1$, let $\phi_1(v)=i_v\omega$; finally, for
$\phi_2:\Lambda^2(V)\rightarrow \Omega^2$, let $\phi_2(v\wedge w)=-(i_vi_w\omega )\omega $.
It is clear that $\phi_0$ and $\phi_1$ are well-defined, while the fact that $\phi_2$ is well-defined follows from
the formula $i_vi_w\omega =-\omega (v\wedge w) \omega$.
Since $\omega$ is nowhere-zero, it follows that $\phi$ is an isomorphism. Observe next that $\phi (\Lambda^2(\jL))=\Lambda^2(d\C [X])$, and since
$X$ is affine we can show that $\Lambda^2(d\C [X])=d\Omega^1 $.
It follows that $\Lambda^2(\dL )/\Lambda^2 (\jL)\cong \Lambda^2 (H^1(X,\C))$, and referring to
the diagram \ref{keycd2} we moreover obtain that $\Lambda^2 (H^1(X,\C))\cong R/\dL \C[X]$, which is (a) and  toghether with 
  part (2) of Proposition \ref{exact} gives also the isomorphism (b). Hence we deduce that 
$\dim \Ker D = \dim \Lambda^2 ({\rm H}^2(X, \C )) -\dim (\C [\omega]\cap \Lambda^2 ({\rm H}^2(X, \C ) ))+\dim \C [\omega$]=
$\dim {\rm span} \{\C [\omega ] , \Lambda^2 ({\rm H}^1 (X, \C))\}$
 and the conclusion of the proof follows from the fact that
 for $n=2$ we have $\Omega^{n-2}/\dL \Omega^{n-2}= \C [X]/\dL \C [X]\cong{\rm  H^2} (X, \C)$ 
\end{proof}

\end{theorem}

\begin{example}
It follows from Theorem  \ref{surface} that $\jL= [\dL , \dL ]$ if $X=\C^2$ or $X=\C^*\times \C^*$.

\end{example}

\begin{example}\label{matias}
If $p(z)$ is a polynomial of degree $d$ with $d$ distinct roots,
then the affine surface $X_p= \{uv =p(z)\}$ is non-singular. Such surfaces have been introduced by Danielewski in order to 
provide counterexamples to the cancellation conjecture.  It is well-known that  $h^1(X)=0$, $h^2(X)=d-1$
and that $X_p$ admits a volume form $\omega$ which is non-exact if and only if $d>1$ (see for example \cite{Le}; $d=1$ then $X_p$ is just the affine plane). 
In particular, for this class of surfaces the derived subalgebra of $\dL $ does not coincide with $\jL$. 
 Let $\lL$ be the Lie algebra generated by locally nilpotent derivations on $X$.
In \cite{Le} the authors show that $\dim \dL/\lL = h^2(X)-1$. Moreover, their calculations show
that $\lL \subset [\dL , \dL ]$, hence $\lL = [\dL , \dL ]$. This equality of infinite-dimensional 
Lie algebras is surprising, since we do not observe any obvious relation 
between $\lL$ and the topology of $X_p$.
\end{example}

\begin{example}
The proof of Theorem \ref{surface} applies to compact Riemann surfaces, while for higher dimensional differentiable manifold the proof of Lichnerowicz is different.
If $X$ is a compact Riemann surface with a given volume form $\omega$, then $ \Lambda^2 ({\rm H}^1 (X, \C))= {\rm H}^2(X, \C)$ hence
$\jL = [\dL , \dL ]$.

\end{example}

\begin{example}
Let $X$ be the complement in $\C^2$ of $n+1$ lines $L_1, ..., L_{n+1}$
containing the origin, and  $\omega$ the volume form obtained by restriction of the volume form on $\C^2$.
It is  proven in \cite{B} that the cohomology ring of $X$ is generated by the closed $1$-forms
$dl_i/l_i$, where  $l_i$ is  the definining polynomial of $L_i$: hence  $ \Lambda^2 ({\rm H}^1 (X, \C))= {\rm H}^2(X, \C)$ and  we conclude that $[\dL , \dL ] =\jL$.

\end{example}

The following theorem, which follows immediately from Proposition \ref{comm} and the commutative diagram \ref{keycd}
will be used to prove the analogue of Theorem \ref{lic} for linear algebraic groups.

\begin{theorem}\label{equal}
If $\Lambda^2 (\dL)=\Lambda^2(V)$, then
$[\dL , \dL ] =\jL$.
\end{theorem}

Next, we recycle a notion
that was introduced in \cite{KK2} in the contest of the Andersen-Lempert theory and density property for affine varieties.

\begin{definition}\label{sc}
Let $X$ be an affine algebraic variety equipped with a volume form.
A pair of  vector fields $\{ v,w \} \subset \dL $ is said to be semi-compatible
if $\Ker v\otimes_{\C} \Ker w$ contains a non-zero ideal of $\C [X]$.
\end{definition}

\begin{theorem}\label{sct}
Let $X$ be an affine algebraic variety equipped with a volume form.
Let $\{ v_k,w_k \}_{k=1}^N$ be a collection of semi-compatible pairs, and for  each $k$  let $I_k \subset \Ker v_k\otimes \Ker w_k$ be a non-trivial ideal as in Definition \ref{sc}.
Suppose that for all $x\in X$ the collection of $\C [X]$-modules
$\{ I_k(v_k\wedge w_k)\}_{k=1}^N$ generate the fiber $\Lambda^2 T_xX$ of the vector bundle $\Lambda^2(TX)$.  Then $[\dL , \dL ] =\jL$.

\begin{proof}
Let $M= \sum_k I_k v_k\wedge w_k$, which is a $\C [X]$-module since $I_k$ are ideals of $\C [X]$.
A regular function $p\in I_k$ can written as a linear combination of product $ab$ of functions
such that $a\in \Ker v_k$ and $b\in \Ker w_k$. Since $av_k\in \dL $ and $bw_k\in \dL$
it follows that $ab(v_k\wedge w_k )=av_k\wedge bw_k\in \Lambda^2 (\dL)$, and therefore $M\subset \Lambda^2(\dL )$.
By hypothesis it follows that $M_x=\Lambda^2(T_xX)$ for all $x\in X$. 
As shown in \cite{KK2}, it follows from Excercise 10, Chapter 2 of \cite{Ha} that $M=\Lambda^2(V)$, and 
  the conclusion of the proof follows therefore from Theorem \ref{equal}.

\end{proof}

\end{theorem}

\begin{example}\label{TA}
Let $X\cong \C^n$ with coordinates $(t_1, t_2,...,t_n)$; then for any $h\neq k$ we have a semi-compatible pair $\{ \frac{\partial}{\partial t_k}, \frac{\partial }{\partial t_h} \}$ satisfying the conditions of Theorem \ref{sct}. Let $X=(\C^*)^n$ equipped with the standard volume form; then for $h\neq k$ we have a
semi-compatible pair $\{ t_k\frac{\partial}{\partial t_k} , t_h\frac{\partial }{\partial t_h} \}$ satyfing the hyphotesis of Theorem \ref{sct}. Hence for such examples $\jL =[\dL , \dL ] $.

\end{example}

Before stating the next lemma, recall that  a vector field on an affine variety $X$ admits a natural extension as a vector field on $X\times Y$, for 
any affine variety $Y$.

\begin{lemma}\label{timescurve}

Let $X$ and $Y$ be equipped respectively with volume forms $\omega$ and $\nu$. 
Let $v,w\in\dL (X)$: if $\{v,w\}$ is semi-compatible pair on $X$, then it is a semi-compatible pair
on $X \times Y$. If $v\in \dL (X)$ and $w\in \dL (Y)$, then $\{v,w\}$ is a semi-compatible pair on $X \times Y$.

\begin{proof}
For a semi-compatible pair  $\{ v ,w\}$ on $X$, let $I\subset\Ker w\otimes \Ker w$ be a non-trivial ideal of $\C [X]$. 
Then $I\otimes \C [Y]\subset\Ker v\otimes \Ker w$ is a non-trivial ideal on $\C [X\times Y]=\C [X]\otimes \C [Y]$.

The second statement follows from the general fact that if $v$ is a vector field on $X$ and $w$ is a vector field on $Y$, then
$\Ker v\otimes \Ker w= \C[X\times Y]$.

\end{proof}

\end{lemma}

We now apply Theorem \ref{equal} to show that the result of Lichnerowicz holds unchanged for all
complex linear algebraic groups of dimension at least $2$.
The proof follows very closely the proof of Theorem 3 from \cite{KK2}. An extra ingredient is necessary however when $G$ is a simple group, where we will need to apply the 
following proposition, whose proof  was suggested by Alexander Premet.

\begin{proposition}\label{semisimple}
Let $G$ be a simple linear algebraic group with Lie algebra $\mathfrak g$. Then
 $\Lambda^2 \mathfrak g$  is spanned by a collection  $\{v_k\wedge w_k \}_{k=1}^N$ such that
for each $k$, $\{ v_k,w_k \}$ is a pair of nilpotent elements corresponding to a  $\mathfrak{sl}_2$-triplet
in $\mathfrak g$.

\begin{proof}

The adjoint action of $G$ on $\mathfrak {g}$ extends naturally to $\Lambda^2 \mathfrak g$.

Consider a pair $\{e , f\}$ of nilpotent element of a  $\mathfrak{sl}_2$-triplet 
where $f$ is  a vector of highest positive root. 
A spanning set can be produced by acting recursively on $e \wedge f$ via
the action of $G$ as follows.

Let  $v$ be a negative root vector and $g=\exp (v)$: then since $[v, e] =0$, it follows that
$Ad_g (e\wedge f ) = e \wedge Ad_g f$.
Let $V$ be the vector space spanned by elements of the form $Ad_{g_1}....Ad_{g_N}f$,
for $g_k= \exp tv_k$ being exponential of negative root vectors $v_k$ and $t\in \C$.
We claim that $V=\mathfrak g$. Since  $\mathfrak g$ is spanned by elements of the form
$[v_1,[v_2,...,[ v_k, f]]...]$, we need to prove that $[v_1,[v_2,..., [v_k, f]]...]\in V$
for any choice of negative root vectors $v_1,...,v_k$. This can be done by induction
on the number of roots applied to $f$.
Since $Ad_{tv}f= f+t[v,f]+ ... +t^N/N! [v,[v,[v,[v ...,f]]]]$,
if we choose $N+1$ values  $t_1, ..., t_{N+1}$ we can
express $[v,f]$ as a linear combination of $Ad{g_k}f$, where $g=\exp ({t_kv})$. Hence $[v,f]\in V$ for all negative root vectors $v$.
Suppose that for any choice of  $k$ negative root vectors $v_1, ... v_k$ it is true that $w=[v_1,[v_2, [...[v_k,f]]]]\in V$,
and let $v_{k+1}$ be another negative root vector.
Then using the same formula as before we can express $[v_{k+1},f]$ as a linear combination of
$Ad_gw$, for $g=\exp (tv_{k+1})$, and since $V$ is invariant under adjoint representation the result follows.
Hence we can find vectors $e\wedge f_1, \cdots e\wedge f_N$,
for $\{e,f_k\}$ a nilpotent pair and
such that $f_1,... f_N$ is a basis of the vector space $\mathfrak g$,
and $f_k=Ad_{g_k}f$, where $g_k$ is a composition of exponentials of negative root vectors.
Next, let $w$ be a positive root vector, and $h=\exp w$.
Then
$Ad_{g_khg_k^{-1}} e\wedge f_k = Ad_{g_k}Ad_h e\wedge f_k$.
Hence, as before, for each fixed $k$
we can find a basis of $\mathfrak g$ of the form $f^k_j$ where
$f^k_j= Ad_{g_k}Ad_{h_j}e$, $h_j=\exp w_j$ is a composition of exponentials of positive root vectors,
and $\{ f_j^k, f_k\}$ is a pair of nilpotent elements in some $\mathfrak {sl}_2$-triplet.

Then it is clear that the set $\{  f_j^k\wedge f_k \}$ is a spanning set of $\Lambda^2 V$
with the desired properties.

\end{proof}

\end{proposition}

\begin{proposition}\label{SL2}
Let $H_1$ and $H_2$ be linear algebraic groups isomorphic to $\C$ acting on a non-singular affine algebraic variety $X$ equipped with a volume form $\omega$
and let $v_1$ and $v_2$ be the vector fields associated respectively with $H_1$ and $H_2$. Suppose that either $[v_1,v_2]=0$ or that 
$v_1$ and $v_2$ are nilpotent elements of the Lie algebra of a group isomorphic to $SL_2$ acting on $X$.
Then $v_1,v_2 \in\dL (X)$ and 
 $\{v_1,v_2\}$ is a compatible pair on $X$. 
Suppose moreover that there exists a finite group
$\Gamma$ acting on $X$ and that $\omega $ is $\Gamma$-invariant. Let $X'=X//\Gamma$ be the quotient of $X$ and 
$\omega'$ be the volume form on $X'$ induced by $\omega$.
If $v_1$ and $v_2$ descends to vector fields $v_1'$ $v_2'$, 
then $v_1',v_2' \in\dL (X')$ and 
 $\{v_1',v_2'\}
$ is a compatible pair on $X'$

\begin{proof}
A vector field corresponding to a non-trivial algebraic  action of $\C$ on an affine algebraic variety equipped with a volume form is zero divergence:
this follows from observing that the ring of regular functions of the general orbits of the action is isomorphic to $\C [t]$, hence has non-trivial units. 
Hence $v_1, v_2\in \dL (X) $ and $v_1',v_2'\in \dL (X')$.
The semicompatibility of the pairs $\{ v_1 , v_2 \}$ and $\{ v_1',v_2'\}$ is proven in Lemma 3.6 of \cite{KK2}.

\end{proof}

\end{proposition}

\begin{theorem}
Let $X$ be a linear algebraic group of dimension at least two. Then $ [ \dL (X), \dL (X) ] = \jL (X)$.

\begin{proof}

If $X$ is one dimensional  then $\dL$ has trivial commutator, hence the theorem is false in one dimension.

As an affine variety $X$ is isomorphic, by Mostov's Theorem, to the product $R\times A$, where $R$ is a maximal closed reductive subgroup of $X$
and $A\cong \C^d$ is the unipotent part of $X$. Moreover, the reductive part $R$ is isomorphic to  $(S\times T)/\Gamma $, for $S$ being a semisimple subgroup of $R$,
$T\cong (\C^*)^t$ and $\Gamma$ a central finite subgroup. Finally, the semisimple part $S$ is isomorphic to the quotient of $S_1\times S_2 ... \times S_N$, for simple groups  $S_1,S_2, ... S_N$
by a central finite subgroup $\Gamma '$.

By Theorem \ref{semisimple} and Proposition \ref{SL2}, for each simple component $S_k$ we have semi-compatible pairs
$\{ v^k_j, w^k_j\}$ such that $v^k_j\wedge w^k_j$  spans $\Lambda^2T_xS_k$ for all $x\in S_k$.  Moreover
for  $h\neq k$ and any $j,l$ the pairs $\{ v^h_j, w^k_l\}$ and $\{ w^h_j, v^k_l\}$ are semi-compatible by Lemma \ref{timescurve}.

Since $\Gamma$ is central, by Lemma \ref{SL2} all the pairs as above descend to semi-compatible pairs on the quotient $S$. 
For convenience of notation, we arrange the semi-compatible pairs, that we constructed on $S$, in the collection $\{\delta_1^j, \delta_2^j\}$.

Next, the vector fields $\sigma_k=t_k\frac{\partial}{\partial t_k}$ are zero divergence vector fields on $T\cong (\C^*)^t$, for each $k=1,...,t$. 
For each $h\neq k$ we have a semi-compatible pair $\{ \sigma_k, \sigma_h\}$, as it follows from Lemma \ref{timescurve}. The same lemma implies that, for all $h\neq k$ ,  $\{ \delta_1^k, \sigma_h \}$ and 
$\{ \delta_2^k, \sigma_h \}$ are semi-compatible pairs. 
Since $\Gamma'$ is central, Lemma \ref{SL2} implies again that  all the  semi-compatiible pairs that we have constructed 
on $S\times T$ descend to semi-compatible pairs on the quotient $(S\times T)/\Gamma$. 
Finally we have vector fields $\mu_k =\frac{\partial }{\partial z_k} $ on $A\cong \C^d$, which give semi-compatible pairs $\{\mu_h, \mu_k\}$ for all 
$k\neq h$.

In conclusion, we constructed the following list of semi-compatible pairs on $X$:
 $\{ \delta^k_1, \sigma_h\}$, $\{\delta_1^k, \mu_j \}$,  $\{ \delta^k_2, \sigma_h\}$, $\{\delta_2^k, \mu_j\}$ ,
$\{\sigma_h, \mu_j \}$, $\{\delta_h, \delta_h\}$ $\{\sigma_j, \sigma_ l\}$, $\{\mu_m,\mu_n \}$.
Since by Theorem \ref{semisimple} $\delta_h\wedge \delta_h$ span $\Lambda^2T_xS$, and $\Gamma $ and $\Gamma '$ are finite groups
it follows that the collection $\{ \delta^k_1\wedge \sigma_h$, $\delta_1^k\wedge \mu_j$,  $\delta^k_2\wedge \sigma_h$, $\delta_2^k\wedge \mu_j$ ,
$\sigma_h\wedge \mu_j $, $\delta_h\wedge \delta_h$ $\sigma_j\wedge \sigma_ l $, $\mu_m\wedge\mu_n \}$
 spans $\Lambda^2T_xG$ for all $x\in G$ and the conclusion 
of the theorem follows from Theorem \ref{equal}.

\end{proof}

\end{theorem}

\vfuzz=2pt

\providecommand{\bysame}{\leavevmode\hboxto3em{\hrulefill}\thinspace}

\end{document}